\def\draftdate{\today}

\documentclass{amsart}
\usepackage{stmaryrd}
\usepackage{xy}
\usepackage{mathrsfs}
\usepackage{amscd,amssymb,color}
\usepackage[colorlinks=true]{hyperref}

\mathchardef\varDelta="7101

\newcommand{\dgcat}{\mathsf{dgcat}}
\newcommand{\dg}{\mathsf{dg}}
\newcommand{\dgHo}{\mathsf{H}^0}
\newcommand{\parf}{\mathsf{parf}}

\renewcommand{\to}{\mathchoice{\longrightarrow}{\rightarrow}{\rightarrow}{\rightarrow}}

\newcommand{\cA}{{\mathcal A}}
\newcommand{\cB}{{\mathcal B}}
\newcommand{\cC}{{\mathcal C}}
\newcommand{\cD}{{\mathcal D}}
\newcommand{\cE}{{\mathcal E}}
\newcommand{\cF}{{\mathcal F}}

\newcommand{\cU}{{\mathcal U}}

\newcommand{\cO}{{\mathcal O}}

\newcommand{\bbK}{I\mspace{-6.mu}K}

\newcommand{\A}{\mathcal{A}}

\newcommand{\bbR}{\mathbb{R}}
\newcommand{\bbZ}{\mathbb{Z}}

\newcommand{\Hmo}{{\mathsf{Hmo}}}

\def\quickop#1{\expandafter\DeclareMathOperator\csname
#1\endcsname{#1}}
\quickop{id}\quickop{Id}\quickop{colim}\quickop{hocolim}\quickop{op}
\quickop{co}\quickop{Ar}\quickop{sing}\quickop{Hom}\quickop{w}\quickop{Ho}
\quickop{ob}\quickop{diag}\quickop{Stab}\quickop{Cat}\quickop{Mot}\quickop{Mod}
\quickop{map}\quickop{Cone}\quickop{End}\quickop{Idem}\quickop{Perf}\quickop{Parf}\quickop{Ind}
\quickop{Gap}\quickop{Shv}\quickop{Spaces}\quickop{cof}\quickop{ex}\quickop{perf}
\quickop{tri}\quickop{Set}\quickop{Fib}\quickop{Nat}\quickop{haut}

\newcommand{\bbL}{\mathbb{L}}

\newcommand{\Trf}{\mathcal{T}r}
\newcommand{\dD}{\mathsf{D}}

\newcommand{\Madd}{\Mot^{\mathsf{add}}_{\dg}}
\newcommand{\Mloc}{\Mot^{\mathsf{loc}}_{\dg}}
\newcommand{\Uadd}{\cU^{\mathsf{add}}_{\dg}}
\newcommand{\Uloc}{\cU^{\mathsf{loc}}_{\dg}}

\newcommand{\ie}{\textsl{i.e.}\ }

\newcommand{\too}{\longrightarrow}

\newcommand{\rep}{\mathrm{rep}} 

\numberwithin{equation}{section}

\newtheorem{theorem}[equation]{Theorem}
\newtheorem*{theorem*}{Theorem}

\newtheorem{lemma}[equation]{Lemma}
\newtheorem{proposition}[equation]{Proposition}
\theoremstyle{definition}

\theoremstyle{remark}

\newtheorem{example}[equation]{Example}

\xyoption{arrow}
\xyoption{matrix}
\xyoption{cmtip}
\SelectTips{cm}{}

\newdir{ >}{{}*!/-5pt/\dir{>}}

\begin{document}

\title[Transfer maps and projection formulas]{Transfer maps and projection formulas}
\author{Gon{\c c}alo~Tabuada}
\address{Departamento de Matem{\'a}tica e CMA, FCT-UNL, Quinta da Torre, 2829-516 Caparica,~Portugal}
\email{tabuada@fct.unl.pt}

\date{\draftdate}


\begin{abstract}
Transfer maps and projection formulas are undoubtedly one of the key tools in the development and computation of (co)homology theories. In this note we develop an unified treatment of transfer maps and projection formulas in the non-commutative setting of dg categories. As an application, we obtain transfer maps and projection formulas in algebraic $K$-theory, cyclic homology, topological cyclic homology, and other scheme invariants.
\end{abstract}

\maketitle
\setcounter{tocdepth}{0}
\section{Transfer maps}
A {\em differential graded (=dg) category}, over a base commutative ring $k$, is a category enriched over complexes of $k$-modules (morphisms sets are complexes)
in such a way that composition fulfills the Leibniz rule\,:
$d(f\circ g)=(df)\circ g+(-1)^{\textrm{deg}(f)}f\circ(dg)$. A {\em dg functor} is a functor which preserves this extra layer of structure; see \S\ref{sec:back}.
Dg categories solve many of the technical problems inherent to triangulated categories and are nowadays widely used in algebraic geometry, representation theory, mathematical physics, $\ldots$; see Keller's ICM adress~\cite{ICM}. We will denote by $\dgcat$ the category of small dg categories.

The purpose of this note is to show that every element of a broad family of functors defined on $\dgcat$ is endowed with an extra degree of ``contravariant functoriality''. Let $\Trf$ be the class of dg functors $F: \cA \rightarrow \cB$ such that for every object $b \in \cB$ the right $\cA$-module $a \mapsto \cB(F(a),b)$ is compact~\cite{Neeman} in the derived category $\cD(\cA)$; see~\S\ref{sub:modules}.

\begin{example}\label{ex:transf}
Natural examples of elements in $\Trf$ are given by morphisms $f: A \to B$ of unital (but not necessary commutative) $k$-algebras such that $B$ is projective of finite type as a right $A$-module; we consider $A$ and $B$ as dg categories with a single object. More generally, we can consider morphisms such that $B$ admits a finite resolution by projective right $A$-modules of finite type. Other examples, of scheme-theoretical nature, are described in Lemma~\ref{lem:schemes}.
\end{example}
\begin{theorem}\label{thm:main}
(i) Let $E:\dgcat \to \dD$ be a functor which sends derived Morita equivalences (see \S\ref{sub:Morita}) to isomorphisms. Then, for every $F \in \Trf$ we have an associated {\em transfer map} $E^{tr}(F): E(\cB) \to E(\cA)$. Moreover this procedure is functorial, \ie $E^{tr}(\id_{\cA}) = \id_{E(\cA)}$ and $E^{tr}(G \circ F) = E^{tr}(F) \circ E^{tr}(G)$.

(ii) Let $\eta: E_1 \Rightarrow E_2$ be a natural transformation between derived Morita invariant functors. Then, $\eta$ is compatible with the transfer maps, \ie for every $F \in \Trf$ we have an equality $\eta_{\cA} \circ E_1^{tr}(F)  = E_2^{tr}(F) \circ \eta_{\cB}$.
\end{theorem}
Connective algebraic $K$-theory ($K$), non-connective algebraic $K$-theory ($\bbK$), Hochschild homology ($HH$), cyclic homology ($HC$), negative cyclic homology ($HN$), periodic cyclic homology ($HP$), topological Hochschild homology ($THH$), and topological cyclic homology ($TC$), are all examples of derived Morita invariant functors; see \cite{BM, Cisinski, Exact, Marco, AGT, Thomason}. Therefore, by Theorem~\ref{thm:main}{\it (i)} we obtain transfer maps\,:
\begin{eqnarray}
&K^{tr}_{\ast}(F): K_{\ast}(\cB) \rightarrow K_{\ast}(\cA) & \bbK^{tr}_{\ast}(F): \bbK_{\ast}(\cB) \rightarrow \bbK_{\ast}(\cA) \label{eq:transf1}\\
&HH^{tr}_{\ast}(F): HH_{\ast}(\cB) \rightarrow HH_{\ast}(\cA) & HC^{tr}_{\ast}(F): HC_{\ast}(\cB) \rightarrow HC_{\ast}(\cA) \label{eq:transf2}\\
& HN^{tr}_{\ast}(F): HN_{\ast}(\cB) \rightarrow HN_{\ast}(\cA) &  HP^{tr}_{\ast}(F): HP_{\ast}(\cB) \rightarrow HP_{\ast}(\cA) \label{eq:transf3}\\
&THH^{tr}_{\ast}(F): THH_{\ast}(\cB) \rightarrow THH_{\ast}(\cA) & TC^{tr}_{\ast}(F): TC_{\ast}(\cB) \rightarrow TC_{\ast}(\cA) \label{eq:transf4}\,.
\end{eqnarray}
When $F$ is as in Example~\ref{ex:transf}, the transfer map $K^{tr}_{\ast}(F)$ \eqref{eq:transf1} coincides with Quillen's original transfer map $f_{\ast}: K_{\ast}(B) \to K_{\ast}(A)$ on algebraic $K$-theory; see \cite[page 103]{Quillen} \cite[\S V 3.3.2]{Weibel}. 
The remaining transfer maps are (to the best of the author's knowledge) new in the literature. It is expected that they will play the same catalytic role that transfer maps did in algebraic $K$-theory.

Now, recall from \cite[\S8.4, \S11.4]{Loday} the construction of the higher Chern characters and the Dennis trace map
\begin{eqnarray}\label{eq:trace}
ch_{\ast, i}: K_{\ast} \Rightarrow HC_{\ast + 2i} & ch^-_{\ast}: K_{\ast} \Rightarrow HN_{\ast} & Dtr_{\ast}: K_{\ast} \Rightarrow THH_{\ast}\,.
\end{eqnarray}
As an important consequence of Theorem~\ref{thm:main}{\it (ii)} the novel transfer maps, introduced in this note, are compatible with the classical natural transformations~\eqref{eq:trace}.

\section{Projection formulas}\label{sec:projection}
The tensor product $-\otimes-$ of $k$-algebras extends naturally to dg categories, giving rise to a symmetric monoidal structure on $\dgcat$; see \S\ref{sub:modules}. In this section we assume that $\dD$ is a symmetric monoidal triangulated category (with unit object ${\bf 1}$), $E:\dgcat \to \dD$ is a symmetric monoidal derived Morita invariant functor, and $F:\cA\to \cB$ is an element in $\Trf$ between dg categories endowed with a ``multiplication''
\begin{eqnarray}\label{eq:mult}
m_{\cA}: \cA \otimes \A \too \cA && m_{\cB}: \cB \otimes \cB \too \cB\,.
\end{eqnarray}
Given an object $P$ in $\dD$, we will denote by $\mathrm{Hom}_{\ast}({\bf 1}, P)$ the $\bbZ$-graded abelian group of morphisms such that $\mathrm{Hom}_n({\bf 1},P):=\mathrm{Hom}(\Sigma^n({\bf 1}),P)$.

Now, let $F: \cA \to \cB$ be an element in $\Trf$ as above. Theorem~\ref{thm:main}{\it (i)} furnish us a transfer map $E^{tr}(F)$, and so by composition with $E(F)$ and $E^{tr}(F)$ we obtain homomorphisms
\begin{eqnarray}
E(F)_{\ast}: \mathrm{Hom}_{\ast}({\bf 1}, E(\cA))& \too & \mathrm{Hom}_{\ast}({\bf 1}, E(\cB)) \label{eq:ind1}\\
E^{tr}(F)_{\ast}: \mathrm{Hom}_{\ast}({\bf 1}, E(\cB))& \too & \mathrm{Hom}_{\ast}({\bf 1}, E(\cA)) \label{eq:ind2}
\end{eqnarray}
of $\bbZ$-graded abelian groups. Since $E$ is symmetric monoidal and the dg categories $\cA$ and $\cB$ are endowed with the dg functors \eqref{eq:mult}, the $\bbZ$-graded abelian groups $\mathrm{Hom}_{\ast}({\bf 1}, E(\cA))$ and $\mathrm{Hom}_{\ast}({\bf 1}, E(\cB))$  naturally become $\bbZ$-graded rings; the multiplicative operation $-\cdot-$ is given by combining the monoidal structure on $\dD$ with the maps $E(m_{\cA})$ and $E(m_{\cB})$.

The category $\dgcat$ carries a Quillen model structure whose weak equivalences are the derived Morita equivalences. The homotopy category obtained will be denoted by $\Hmo$; see \S\ref{sub:Morita} for a description of its $\mathrm{Hom}$-sets in terms of bimodules. The tensor product of dg categories can be derived $-\otimes^{\bbL}-$, making $\Hmo$ into a symmetric monoidal category; with unit object the dg category $\underline{k}$ with one object and $k$ as the dg algebra of endomorphisms. Thanks to Theorem~\ref{thm:main}{\it (i)} the functor $\dgcat \to \Hmo$ is endowed with transfer maps. We will denote by $[X_F]$ the image of $F$ in $\Hmo$ and by $[{}_{F}X]$ the associated transfer map. Under the above assumptions we have the following result.
\begin{proposition}\label{prop:proj}
Let $F$ be an element in $\Trf$ making the following diagrams commutative
\begin{equation}\label{eq:diagram}
\xymatrix@C=2em@R=3em{
\cA \otimes^{\bbL} \cA \ar[d]_{m_{\cA}} \ar[rr]^-{[X_F] \otimes^{\bbL} [X_F]} && \cB \otimes^{\bbL} \cB \ar[d]^{m_{\cB}} & \cB \otimes^{\bbL}\cA \ar[d]_{\id \otimes^{\bbL}[X_F]} \ar[rr]^{[{}_FX]\otimes^{\bbL}\id} && \cA \otimes^{\bbL} \cA \ar[r]^-{m_{\cA}} & \cA \\
\cA \ar[rr]_{[X_F]} && \cB & \cB \otimes^{\bbL}\cB \ar[rrr]_{m_{\cB}} &&& \cB \ar[u]_{[{}_FX]}   \,.
}
\end{equation}   
Then, for every $x \in \mathrm{Hom}_{\ast}({\bf 1},E(\cB))$ and $y \in \mathrm{Hom}_{\ast}({\bf 1},E(\cA))$, the following {\em projection formula}
\begin{equation}\label{eq:proj}
E^{tr}(F) \circ (x \cdot (E(F) \circ y)) = (E^{tr}(F) \circ x) \cdot y 
\end{equation}
holds in $\mathrm{Hom}_{\ast}({\bf 1}, E(\cA))$.
\end{proposition}
Note that the commutativity of the left-hand side diagram in \eqref{eq:diagram} implies that \eqref{eq:ind1} is a homomorphism of $\bbZ$-graded rings. In this case, the $\bbZ$-graded ring $\mathrm{Hom}_{\ast}({\bf 1},E(\cB))$ becomes a $\mathrm{Hom}_{\ast}({\bf 1}, E(\cA))$-algebra and so the projection formula \eqref{eq:proj} entails that the homomorphism \eqref{eq:ind2} is $\mathrm{Hom}_{\ast}({\bf 1}, E(\cA))$-linear. 
\begin{example}\label{ex:example2}
Natural examples of elements in $\Trf$ which satisfy the conditions of Proposition~\ref{prop:proj} are given by the morphisms $f:A \to B$ between {\em commutative} $k$-algebras presented in Example~\ref{ex:transf}. In this commutative case the multiplication operations in $A$ and $B$ are $k$-algebra homomorphisms, and so we obtain dg functors as in \eqref{eq:mult}. The commutativity of the left-hand side diagram in \eqref{eq:diagram} is clear. The commutativity of the right-hand side diagram in \eqref{eq:diagram} follows from the natural isomorphism $B\otimes_BA \otimes_A B \simeq B \otimes_A A$ of right $A$-modules. Scheme-theoretical examples are described in Proposition~\ref{prop:schemes}.
\end{example}  
We now give examples of symmetric monoidal derived Morita invariant functors.
\begin{example}[Algebraic $K$-theory]\label{ex:1}
Recall from \cite[\S10, \S15]{Duke} the construction of the universal additive and localizing invariant of dg categories 
\begin{eqnarray*}
\Uadd: \dgcat \too \Madd(e) && \Uloc: \dgcat \to \Mloc(e)\,.
\end{eqnarray*}
Here, $\Madd(e)$ and $\Mloc(e)$ are the categories of {\em non-commutative motives} in the additive and localizing settings; see \cite{CT, CT1, Duke} for details. These triangulated categories carry symmetric monoidal structures with unit objects $\Uadd(\underline{k})$ and $\Uloc(\underline{k})$. Thanks to \cite[Theorem~15.10]{Duke} and \cite[Theorem~7.16]{CT} we have, for every small dg category $\cC$, natural isomorphisms of $\bbZ$-graded abelian groups
\begin{eqnarray*}
\mathrm{Hom}_{\ast}(\Uadd(\underline{k}), \Uadd(\cC))\simeq K_{\ast}(\cC) && \mathrm{Hom}_{\ast}(\Uloc(\underline{k}), \Uloc(\cC))\simeq \bbK_{\ast}(\cC)\,.
\end{eqnarray*}
Therefore, it follows from Proposition~\ref{prop:proj} with $E=\Uadd$, that $x \in K_{\ast}(\cB)$, $y \in K_{\ast}(\cA)$, the map \eqref{eq:ind2} identifies with the transfer map $K^{tr}_{\ast}(F)$ \eqref{eq:transf1}, and the projection formula \eqref{eq:proj} holds in $K_{\ast}(\cA)$. Similarly for $E=\Uloc$, with $K$ replaced by $\bbK$.

When $F$ is an in Example~\ref{ex:example2}, the multiplicative operation on $K_{\ast}(A)$ and $K_{\ast}(B)$ is the one described in \cite[\S11.2.16]{Loday}. In this case we recover Quillen's original projection formula $f_{\ast}(x \cdot f^{\ast}(y)) = f_{\ast}(x) \cdot y$ on algebraic $K$-theory; see \cite[page 103]{Quillen} \cite[\S V 3.5.3]{Weibel}.  
In the non-connective algebraic $K$-theory case, and in all the following examples, we obtain analogous projection formulas which are new in the literature.
\end{example}
\begin{example}[Hochschild homology]
Recall from example \cite[Example~7.9]{CT1} the construction of the Hochschild homology functor $HH : \dgcat \to \cD(k)$. For every small dg category $\cC$ we have a natural isomorphism of $\bbZ$-graded $k$-modules $\mathrm{Hom}_{\ast}(k, HH(\cC)) \simeq HH_{\ast}(\cC)$. Therefore, it follows from Proposition~\ref{prop:proj} with $E=HH$, that $x \in HH_{\ast}(\cB)$, $y \in HH_{\ast}(\cA)$, the map \eqref{eq:ind2} identifies with the transfer map $HH^{tr}_{\ast}(F)$ \eqref{eq:transf2}, and the projection formula \eqref{eq:proj} holds in $HH_{\ast}(\cA)$. Moreover, when $F$ is as in Example~\ref{ex:example2}, the multiplicative operation on $HH_{\ast}(A)$ and $HH_{\ast}(B)$ is given by the shuffle product; see \cite[Corollary~4.2.7]{Loday}.
\end{example}
\begin{example}[Negative cyclic homology]
Recall from \cite[Example~7.10]{CT1} the construction of the mixed complex functor $C: \dgcat \to \cD(\Lambda)$. Here, $\Lambda$ is the $k$-algebra $k[\epsilon]/\epsilon^2$ with $\epsilon$ of degree $-1$ and $d(\epsilon)=0$. The derived category $\cD(\Lambda)$ carries a symmetric monoidal structure defined on the underlying complexes, with unit object $k$. For every small dg category $\cC$ we have a natural isomorphism of $\bbZ$-graded $k$-modules $\mathrm{Hom}_{\ast}(k, C(\cC)) \simeq HN_{\ast}(\cC)$. Therefore, it follows from Proposition~\ref{prop:proj} with $E=C$, that $x \in HN_{\ast}(\cB)$, $y \in HN_{\ast}(\cA)$, the map \eqref{eq:ind2} identifies with the transfer map $HN^{tr}_{\ast}(F)$ \eqref{eq:transf3}, and the projection formula \eqref{eq:proj} holds in $HN_{\ast}(\cA)$. Moreover, when $F$ is as in Example~\ref{ex:example2}, the multiplicative operation on $HN_{\ast}(A)$ and $HN_{\ast}(B)$ is the one described in \cite[\S5.1.13]{Loday} in terms of shuffle products.
\end{example}
\begin{example}[Periodic cyclic homology]\label{ex:4}
Recall from \cite[Example~7.11]{CT1} the construction of the periodic complex functor $P: \cD(\Lambda) \to \cD(k[u]\text{-}\mathrm{Comod})$. Here, $k[u]\text{-}\mathrm{Comod}$ is the category of $k[u]$-comodules, where $k[u]$ is the cocomutative Hopf algebra of polynomials in one variable $u$ of degree $2$. The derived category $\cD(k[u]\text{-}\mathrm{Comod})$ carries a symmetric monoidal structure given by the cotensor product of comodules, with unit object $k[u]$. For every small dg category $\cC$ we have a natural isomorphism of $\bbZ$-graded $k$-modules $\mathrm{Hom}_{\ast}(k, (P\circ C)(\cC)) \simeq HP_{\ast}(\cC)$. Therefore, it follows from Proposition~\ref{prop:proj} with $E=P\circ C$, that $x \in HP_{\ast}(\cB)$, $y \in HP_{\ast}(\cA)$, the map \eqref{eq:ind2} identifies with the transfer map $HP^{tr}_{\ast}(F)$ \eqref{eq:transf3}, and the projection formula \eqref{eq:proj} holds in $HP_{\ast}(\cA)$. Moreover, when $F$ is as in Example~\ref{ex:example2}, the multiplicative operation on $HP_{\ast}(A)$ and $HP_{\ast}(B)$ is the one described in \cite[\S5.1.13]{Loday} in terms of shuffle products.
\end{example}
\section{Schemes}\label{sec:schemes}
In this section we assume that our base commutative ring $k$ is a field. By a {\em scheme} we mean a quasi-compact and quasi-separated $k$-scheme; see \cite{SGA6}.

Given a scheme $X$, we have a Grothendieck category $\mathrm{QCoh}(X)$ of quasi-coherent sheaves of $\cO_X$-modules and an associated dg category $\cC_{\dg}(\mathrm{QCoh}(X))$ of unbounded cochain complexes on $\mathrm{QCoh}(X)$. Let $\cC_{\parf}^{\dg}(X)$ be the full small dg subcategory of $\cC_{\dg}(\mathrm{QCoh}(X))$ consisting of those perfect complexes of $\cO_X$-modules (see \cite[\S2]{Thomason}) which are bounded above, degreewise flat, and whose stalks have cardinality at most the cardinality of $k$. We will denote by $\cD_{\parf}^{\dg}(X)$ the Drinfeld's DG quotient~\cite{Drinfeld} (see also \cite{JAlg}) of $\cC_{\parf}^{\dg}(X)$ by its full dg subcategory of acyclic complexes; note that since we are working over a field, the homotopical flatness conditions~\cite[\S3.3]{Drinfeld} are automatically fulfilled. The small dg category $\cD_{\parf}^{\dg}(X)$ enhances the derived category $\cD_{\parf}(X)$ of perfect complexes of $\cO_X$-modules in the sense that we have a natural equivalence of triangulated categories $\dgHo(\cD_{\parf}^{\dg}(X))\simeq \cD_{\parf}(X)$; see \cite{Orlov} for the uniqueness of this dg-enhancement. Given a morphism of schemes $f: X \to Y$ we have a dg-enhancement $f^{\ast}: \cD_{\parf}^{\dg}(Y) \to \cD_{\parf}^{\dg}(X)$ of the derived inverse image functor $\bbL f^{\ast}: \cD_{\parf}(Y) \to \cD_{\parf}(X)$.
\begin{lemma}\label{lem:schemes}
Let $f:X \to Y$ be a perfect projective morphism, a flat proper morphism, or a proper morphism of finite $Tor$-dimension between noetherian schemes; see \cite{SGA6}\cite[\S3]{Thomason}. Then, the dg functor $f^{\ast}: \cD^{\dg}_{\parf}(Y) \to \cD^{\dg}_{\parf}(X)$ belongs to $\Trf$.
\end{lemma}
Thanks to the work of Thomason-Trobaugh~\cite[\S3]{Thomason}, Schlichting~\cite[\S8 Theorem~5]{Marco}, Keller~\cite[\S5.2]{Exact1}, and Blumberg-Mandell~\cite[Theorem~1.3]{BM}, the algebraic $K$-theory, the cyclic homology (and all its variants), and even the topological cyclic homology of a scheme $X$ can be recovered from its associated dg category $\cD_{\parf}^{\dg}(X)$. Therefore, by combining Lemma~\ref{lem:schemes} with \eqref{eq:transf1}-\eqref{eq:transf4}, we obtain transfer maps\,:
$$
\begin{array}{lcl}
f_{\ast}: K_{\ast}(X) \to K_{\ast}(Y) && f_{\ast}: \bbK_{\ast}(X) \to \bbK_{\ast}(Y) \\
f_{\ast}: HH_{\ast}(X) \to HH_{\ast}(Y) && f_{\ast}: HC_{\ast}(X) \to HC_{\ast}(Y) \\
f_{\ast}: HN_{\ast}(X) \to HN_{\ast}(Y) &&  f_{\ast}: HP_{\ast}(X) \to HP_{\ast}(Y)\\
f_{\ast}: THH_{\ast}(X) \to THH_{\ast}(Y) && f_{\ast}: TC_{\ast}(X) \to TC_{\ast}(Y)\,.
\end{array}
$$
The transfer map $f_{\ast}:K_{\ast}(X) \to K_{\ast}(Y)$ coincides with the one originally defined by Quillen \cite[\S7 2.7]{Quillen} and Thomason-Trobaugh~\cite[\S3.16.4-3.16.6]{Thomason}. The remaining transfer maps are new in the literature. Moreover, thanks to Theorem~\ref{thm:main} and \eqref{eq:trace}, these novel transfer maps are compatible with the higher Chern characters and the Dennis trace map.
\begin{proposition}\label{prop:schemes}
Let $f:X \to Y$ be a morphism of schemes as in Lemma~\ref{lem:schemes}. Then, the dg functor $f^{\ast}: \cD^{\dg}_{\parf}(Y) \to \cD^{\dg}_{\parf}(X)$ satisfies the conditions of Proposition~\ref{prop:proj}.
\end{proposition}
By combining Proposition~\ref{prop:schemes} with Examples~\ref{ex:1}-\ref{ex:4}, we obtain the following projection formula
\begin{equation*}\label{eq:proj1}
f_{\ast}(x\cdot f^{\ast}(y)) = f_{\ast}(x) \cdot y\,,
\end{equation*}
which holds in $K_{\ast}(Y)$, $\bbK_{\ast}(Y)$, $HH_{\ast}(Y)$, $HN_{\ast}(Y)$ and $HP_{\ast}(Y)$. In the connective  algebraic $K$-theory case we recover the projection formula originally developed by Quillen \cite[Proposition~2.10]{Quillen} and Thomason-Trobaugh \cite[Proposition~3.17]{Thomason}. In all the remaining cases no projection formulas have been available so far. Hence, Proposition~\ref{prop:schemes} and Lemma~\ref{lem:schemes} shed novel light over the (contravariant) functoriality of these classical invariants of schemes. It is excepted that these new results play the same catalytic role that the projection formula and the transfer maps did in algebraic $K$-theory.

\section{Background on dg categories}\label{sec:back}
Let $\cC(k)$ be the category of (unbounded) cochain complexes of $k$-modules. A {\em differential graded (=dg) category} (over the base ring $k$) is a $\cC(k)$-category and a {\em dg functor} is a $\cC(k)$-functor; see \cite[Definitions~6.2.1-6.2.3]{Borceaux}. For a survey article on dg categories we invite the reader to consult Keller's ICM adress~\cite{ICM}. The category of small dg categories is denoted by $\dgcat$.

\subsubsection{(Bi)modules}\label{sub:modules}
Let $\cA$ be a small dg category. The category $\dgHo(\cA)$ has the same objects as $\cA$ and morphisms given by $\dgHo(\cA)(x,y):=\textrm{H}^0(\cA(x,y))$. The {\em opposite} dg category $\mathcal{A}^{\op}$ of $\cA$ has the same objects as $\mathcal{A}$ and complexes of morphisms given by $\mathcal{A}^{\op}(a_1,a_2):=\mathcal{A}(a_2,a_1)$. A {\em right $\cA$-module} (or simply a $\cA$-module) is a dg functor
$\cA^{\op} \to \cC_{\dg}(k)$ with values in the dg category $\cC_{\dg}(k)$ of complexes of $k$-modules. We denote by $\cC(\cA)$ the category of $\cA$-modules; see \cite[\S 2.3]{ICM}. The {\em derived category $\cD(\cA)$ of $\cA$} is the localization of $\cC(\cA)$ with respect to the class of objectwise quasi-isomorphisms. As explained in \cite[\S 3.1]{ICM} the differential graded structure of $\cC_{\dg}(k)$ makes $\cC(\cA)$ naturally into a dg category $\cC_\dg(\cA)$.

Now, let $\cA$ and $\cB$ be two small dg categories. Their {\em tensor product} $\cA\otimes \cB$ is defined as follows\,: the set of objects is the cartesian product of the sets of objects of $\cA$ and $\cB$ and for every two  objects $(a_1,b_1)$ and $(a_2,b_2)$ in $\cA\otimes \cB$ we have $\cA\otimes\cB((a_1,b_1),(a_2,b_2)):=\cA(a_1,a_2) \otimes \cB(b_1,b_2)$. The tensor product of dg categories gives rise to a symmetric monoidal structure on $\dgcat$, with unit object the dg category $\underline{k}$ with one object and $k$ as the dg algebra of endomorphisms. Finally, by a {\em $\cA\text{-}\cB$-bimodule} we mean a dg functor $\cA \otimes\cB^{\op} \to \cC_{\dg}(k)$.

\subsubsection{Derived Morita equivalences}\label{sub:Morita}
A dg functor $F: \cA \to \cB$ is a called a {\em derived Morita equivalence}
if its derived extension of scalars functor $\bbL F_!: \cD(\cA) \stackrel{\sim}{\to}~\cD(\cB)$ is an equivalence of categories; see \cite[\S5]{IMRN}. Thanks to \cite[Theorem~5.3]{IMRN} the category $\dgcat$ carries a Quillen model structure~\cite{Quillen1} whose weak equivalences are the derived Morita equivalences. We denote by $\Hmo$ the homotopy category hence obtained. The tensor product of dg categories can be derived into a bifunctor $-\otimes^{\bbL}-$, making $\Hmo$ into a symmetric monoidal category; see \cite[\S4]{Toen}.

Now, let $\cA$ and $\cB$ be two small dg categories. We denote by $\rep(\cA, \cB)$ the full sucategory of $\cD(\cA^{\op} \otimes^{\bbL} \cB)$ consisting of those $\cA\text{-}\cB$-bimodules $X$ such that for every object $a\in \cA$ the right $\cB$-module $X(-,a)$ is a compact object~\cite{Neeman} in the derived category $\cD(\cB)$. Thanks to \cite[Corollary]{Toen} and \cite[Remark~5.11]{IMRN} we have a natural bijection $\mathrm{Hom}_{\Hmo}(\cA,\cB) \simeq \mathrm{Iso}\,\rep(\cA,\cB)$, where $\mathrm{Iso}$ denotes the set of isomorphism classes. Moreover, given small dg categories $\cA$, $\cB$ and $\cC$, the composition operation in $\Hmo$ corresponds to the derived tensor product of bimodules
\begin{eqnarray*} 
\mathrm{Iso}\,\rep(\cA,\cB) \times \mathrm{Iso}\,\rep(\cB,\cC) \too \mathrm{Iso}\,\rep(\cA,\cC) && ([X],[Y]) \mapsto [X \otimes^{\bbL}_{\cB} Y]\,.
\end{eqnarray*}
The localization functor $\dgcat \to \Hmo$ is the identity on objects and sends a dg functor $F:\cA \to \cB$ to the isomorphism class of the $\cA\text{-}\cB$-bimodule 
\begin{eqnarray*}
X_F: \cA \otimes^{\bbL}\cB^{\op} \too \cC_{\dg}(k) && (a,b) \mapsto \cB(b,F(a))\,.
\end{eqnarray*}
\section{Proofs}
\begin{proof}{(Theorem~\ref{thm:main})}\label{proof:1}
The homotopy category $\Hmo$ (see \S\ref{sub:Morita}) is the localization of $\dgcat$ with respect to the class of derived Morita equivalences. Therefore, since the functor $E$ sends derived Morita equivalences to isomorphisms, we obtain a well-defined functor $\overline{E}$ making the diagram
\begin{equation}\label{eq:diag}
\xymatrix{
\dgcat \ar[d] \ar[r]^-E & \dD \\
\Hmo \ar[ur]_{\overline{E}} &
}
\end{equation} 
commute. Now, let $F:\cA \to \cB$ be an arbitrary dg functor. Associated to $F$ we have the following $\cB\text{-}\cA$-bimodule
\begin{eqnarray}\label{eq:bimodule}
{}_FX: \cB \otimes^{\bbL}\cA^{\op} \too \cC_{\dg}(k) && (b,a) \mapsto \cB(F(a),b)\,.
\end{eqnarray}
Note that $F$ belongs to $\Trf$ if and only if the bimodule ${}_FX$ belongs to $\rep(\cB,\cA)$. In particular, if $F$ belongs to $\Trf$ we obtain a well-defined morphism $[{}_FX]$ in $\mathrm{Hom}_{\Hmo}(\cB,\cA)=\mathrm{Iso}\, \rep(\cB,\cA)$. The transfer map $E^{tr}(F)$, associated to $F$, is then the morphism in $\dD$ given by $\overline{E}([{}_FX])$. If $G: \cB \to \cC$ is another dg functor, we have a natural isomorphism ${}_{(G\circ F)}X \simeq {}_FX \otimes^{\bbL}_{\cB}{}_GX$ of $\cC\text{-}\cA$-bimodules. In particular, if $F$ and $G$ belong to $\Trf$ we conclude that ${}_{(G\circ F)}X$ belongs to $\rep(\cC, \cA)$. We obtain then the equality $[{}_{(G \circ F)}X]=[{}_FX] \circ [{}_GX]$ in $\Hmo$, which implies the equality $E^{tr}(G \circ F)=E^{tr}(F) \circ E^{tr}(G)$ in $\dD$. This shows item {\it (i)}. 

Let us now show item {\it (ii)}.  Since the functors $E_1$ and $E_2$ are derived Morita invariant, the natural transformation $\eta$ between $E_1$ and $E_2$ is also a natural transformation between the induced functors $\overline{E_1}$ and $\overline{E_2}$; see the above diagram~\eqref{eq:diag}. Therefore, the equalities
\begin{eqnarray*}
E_1^{tr}(F)=\overline{E_1}([{}_FX]) && E_2^{tr}(F)=\overline{E_2}([{}_FX]) 
\end{eqnarray*}
allow us to conclude that $\eta$ is compatible with the transfer maps, \ie for every $F \in \Trf$ we have an equality $\eta_{\cA} \circ E_1^{tr}(F)  = E_2^{tr}(F) \circ \eta_{\cB}$.
\end{proof}
\begin{proof}{(Proposition~\ref{prop:proj})}
Since the functor $E$ sends derived Morita equivalences to isomorphisms, we obtain a well-defined functor $\overline{E}$ as in the above diagram~\eqref{eq:diag}. Moreover, since by hypothesis $E$ is symmetric monoidal so it is $\overline{E}$. Therefore, when we apply the functor $\overline{E}$ to the right-hand side diagram in \eqref{eq:diagram} and use the equalities
\begin{eqnarray*}
E(F)=\overline{E}([X_F]) && E^{tr}(F)=\overline{E}([{}_FX]) \,,
\end{eqnarray*} 
we obtain a commutative diagram in $\dD$\,:
\begin{equation}\label{eq:diag2}
\xymatrix{
E(\cB) \otimes E(\cA) \ar[d]_{\id \otimes E(F)} \ar[rr]^-{E^{tr}(F) \otimes \id} && E(\cA) \otimes E(\cA) \ar[rr]^-{E(m_{\cA})} && E(\cA) \\
E(\cB) \otimes E(\cB) \ar[rrrr]_{E(m_{\cB})} &&&& E(\cB) \ar[u]_{E^{tr}(F)} \,.
}
\end{equation}
Now, let $x \in \mathrm{Hom}_n({\bf 1}, E(\cB))$ and $y \in \mathrm{Hom}_m({\bf 1}, E(\cA))$. Using the symmetric monoidal structure on $\dD$ we obtain an element $x \otimes y$ in $\mathrm{Hom}_{n+m}({\bf 1}, E(\cB)\otimes E(\cA))$. Thanks to the above diagram \eqref{eq:diag2} we have the following equality
\begin{equation}\label{eq:equality}
E^{tr}(F) \circ E(m_{\cB}) \circ (\id \otimes E(F)) \circ (x \otimes y) = E(m_{\cA}) \circ (E^{tr}(F) \otimes \id) \circ (x \otimes y)\,.
\end{equation}
Recall that the multiplicative operation $-\cdot -$ on $\mathrm{Hom}_{\ast}({\bf 1}, E(\cB))$ and $\mathrm{Hom}_{\ast}({\bf 1}, E(\cA))$ is given by combining the monoidal structure on $\dD$ with the maps $E(m_{\cA})$ and $E(m_{\cB})$. Therefore, we observe that the above equality \eqref{eq:equality} coincides with the projection formula
\begin{equation*}
E^{tr}(F) \circ (x \cdot (E(F) \circ y)) = (E^{tr}(F) \circ x) \cdot y \,.
\end{equation*}
This achieves the proof.
\end{proof}
\begin{proof}{(Lemma~\ref{lem:schemes})}
Recall form the proof of Theorem~\ref{thm:main} that the dg functor $f^{\ast}:\cD_{\parf}^{\dg}(Y) \to \cD_{\parf}^{\dg}(X)$ belongs to $\Trf$ if and only if the bimodule ${}_{f^{\ast}}X$~\eqref{eq:bimodule} belongs to $\rep(\cD_{\parf}^{\dg}(X), \cD_{\parf}^{\dg}(Y))$. By tensoring with the bimodule ${}_{f^{\ast}}X$ we obtain a dg functor (see \cite[\S3.8]{ICM})
\begin{equation}\label{eq:tensoring}
-\otimes{}_{f^{\ast}}X : \cC_{\dg} (\cD_{\parf}^{\dg}(X)) \too \cC_{\dg} (\cD_{\parf}^{\dg}(Y))\,.
\end{equation}
If $f$ is as in Lemma~\ref{lem:schemes}, then the derived direct image functor
\begin{equation}\label{eq:direct}
\bbR f_{\ast}: \cD(X) \too \cD(Y)
\end{equation}
preserves perfect complexes (see \cite[\S3.16]{Thomason}) and so the left derived functor of \eqref{eq:tensoring} identifies with \eqref{eq:direct}. Since for every perfect complex $\cE^{\cdot}$ of $\cO_X$-modules we have a natural isomorphism $\cE^{\cdot} \otimes {}_{f^{\ast}}X \simeq \bbR f_{\ast}(\cE^{\cdot})$ in the derived category $\cD_{\parf}(Y)$, we conclude that the bimodule ${}_{f^{\ast}}X$ belongs to $\rep(\cD_{\parf}^{\dg}(X), \cD_{\parf}^{\dg}(Y))$. This achieves the proof.
\end{proof}
\begin{proof}{(Proposition~\ref{prop:schemes})}
The derived inverse image functor
$$\bbL f^{\ast} : \cD_{\parf}(Y) \too \cD_{\parf}(X)$$
is symmetric monoidal with respect to the derived tensor products $-\otimes^{\bbL}_{\cO_Y}-$ and $-\otimes^{\bbL}_{\cO_X}-$; see \cite[\S3.15]{Thomason}. By lifting these derived tensor products to $\cD_{\parf}^{\dg}(Y)$ and $\cD_{\parf}^{\dg}(X)$ we observe that the left-hand side diagram in \eqref{eq:diagram} (with $\cA=\cD_{\parf}^{\dg}(Y)$, $\cB=\cD_{\parf}^{\dg}(X)$, and $F=f^{\ast}$) commutes. Now, let $\cE^{\cdot}$ be a perfect complex of $\cO_X$-modules and $\cF^{\cdot}$ a perfect complex of $\cO_Y$-modules. Thanks to \cite[\S III 3.7]{SGA6} we have a canonical isomorphism
\begin{equation}\label{eq:can}
\bbR f_{\ast} (\cE^{\cdot}) \otimes^{\bbL}_{\cO_Y} \cF^{\cdot} \simeq \bbR f_{\ast} (\cE^{\cdot} \otimes^{\bbL}_{\cO_X} \bbL f^{\ast}(\cF^{\cdot}))
\end{equation}
in the derived category $\cD_{\parf}(Y)$. The left-hand side of \eqref{eq:can} identifies with the evaluation of the bimodule $m_{\cA} \circ [{}_FX] \otimes^{\bbL} \id$ at $\cE^{\cdot} \otimes \cF^{\cdot}$. The right-hand side of \eqref{eq:can} identifies with the evaluation of the bimodule $[{}_FX]\circ m_{\cB} \circ \id \otimes^{\bbL} [X_F]$ at $\cE^{\cdot}\otimes \cF^{\cdot}$. Therefore, the right-hand side diagram in \eqref{eq:diagram} commutes and so the proof is finished.
\end{proof}
\medbreak\noindent\textbf{Acknowledgments:} The author is very grateful to Andrei Suslin for bringing the importance of developing a unified treatment of transfer maps (and associated projection formulas) to his attention. The author would also like to thank the Department of Mathematics at UIC (where this work was carried out) for its hospitality and the Midwest Topology Network for financial support.


\begin{thebibliography}{00}

\bibitem{BM} A.~Blumberg and M.~Mandell, {\em Localization theorems in topological
Hochschild homology and topological cyclic homology}. Available at arXiv:$0802.3938$.    

\bibitem{Borceaux} F.~Borceux, {\em Handbook of categorical algebra. 2}. Encyclopedia of Mathematics and its Applications {\bf 51} (1994). Cambridge Univ. Press.
 
\bibitem{Cisinski} D.-C.~Cisinski, {\em Invariance de la $K$-th{\'e}orie par {\'e}quivalences d{\'e}riv{\'e}es}. Available at  {\tt www-math.univ-paris13.fr/cisinski/publications.html}. To appear in J. of $K$-theory. 
 
\bibitem{CT} D.-C.~Cisinski and G.~Tabuada, {\em Non-connective $K$-theory via universal invariants}.
Available at arXiv:$0903.3717\textrm{v}2$. 
 
\bibitem{CT1} \bysame, {\em Symmetric monoidal structure on Non-commutative motives}.
Available at arXiv:$1001.0228\textrm{v}2$.

\bibitem{Drinfeld} V.~Drinfeld, {\em DG quotients of DG categories}.
J. Algebra {\bf 272} (2004), 643--691.

\bibitem{SGA6} A.~Grothendieck, {\em Th{\'e}orie des intersections et th{\'e}or{\`e}me de Riemann-Roch}. S{\'e}minaire de G{\'e}om{\'e}trie Alg{\'e}brique du Bois-Marie 1966--1967 (SGA 6). Lecture Notes in Mathematics {\bf 225} (1971).
 
\bibitem{ICM} B.~Keller, {\em On differential graded categories}. International Congress of Mathematicians (Madrid), Vol.~{\bf II} (2006), 151--190. Eur.~Math.~Soc., Z{\"u}rich.

\bibitem{Exact} \bysame, {\em On the cyclic homology of exact categories}. J.~Pure~Appl.~Alg. {\bf 136}(1) (1999), 1--56.
    
\bibitem{Exact1} \bysame, {\em On the cyclic homology of ringed spaces and schemes}. Doc. Math. {\bf 3} (1998), 231--259.

\bibitem{Loday} J.-L.~Loday, {\em Cyclic homology}. Grundlehren der Mathematischen Wissenschaften {\bf 301} (1992). Springer-Verlag, Berlin.

\bibitem{Orlov} V.~Lunts and D.~Orlov {\em Uniqueness of enhancement for triangulated categories}. J. Amer. Math. Soc. {\bf 23} (2010), 853--908. 

\bibitem{Neeman} A.~Neeman, {\em Triangulated categories}. Ann.~Math.~Studies {\bf 148} (2001). Princeton Univ. Press.

\bibitem{Quillen} D.~Quillen, {\em Higher algebraic $K$-theory I}. Lecture Notes in Mathematics  {\bf 341} (1973), 85--147.

\bibitem{Quillen1} \bysame, {\em Homotopical algebra}. Lecture Notes in Mathematics {\bf 43} (1967).

\bibitem{Marco} M.~Schlichting, {\em Negative $\mbox{K}$-theory of derived categories}. Math.~Z. {\bf 253} (2006), no.~1, 97--134.

\bibitem{Duke} G.~Tabuada, {\em Higher $K$-theory via universal invariants}.
Duke Math. J. {\bf 145} (2008), no.~1, 121--206.

\bibitem{IMRN} \bysame, {\em Invariants additifs de dg-cat{\'e}gories}.
Int.~Math.~Res.~Not. {\bf 53} (2005), 3309--3339.

\bibitem{JAlg} \bysame, {\em On Drinfeld's DG quotient}. J. Algebra {\bf 323} (2010), 1226--1240.

\bibitem{AGT} \bysame, {\em Generalized spectral categories, topological
Hochschild homology, and trace maps}. Algebraic and Geometric Topology {\bf 10} (2010), 137--213.

\bibitem{Thomason} R.~W.~Thomason and T.~Trobaugh,
{\em Higher algebraic $K$-theory of schemes and of derived 
categories}. Grothendieck Festschrift, Volume III. Volume {\bf 88} of
Progress in Math., 247--436. Birkhauser, Boston, Bassel, Berlin, 1990. 

\bibitem{Toen} B.~To{\"e}n, {\em The homotopy theory of dg-categories and
    derived Morita theory}. Invent. Math. {\bf 167} (2007), no.~3, 615--667.

\bibitem{Weibel} C.~Weibel, {\em The $K$-book: an introduction to algebraic $K$-theory}. A graduate textbook in progress. Available at {\tt http://www.math.rutgers.edu/weibel/Kbook.html}.


\end{thebibliography}
\end{document}